\documentclass[12pt]{article}

\usepackage{lmodern}

\usepackage{scrextend}

\usepackage[left=25mm, right=25mm, top=20mm, bottom=20mm]{geometry}

\usepackage{amsmath,amssymb,amsfonts,setspace,array}
\usepackage{hyperref}
\usepackage{amsthm,graphicx,amstext,mathrsfs,float,lineno}
\theoremstyle{definition}
\newtheorem{theorem}{Theorem}[section] 
\newtheorem{example}[theorem]{Example}
\newtheorem{definition}[theorem]{Definition}
\newtheorem{proposition}[theorem]{Proposition}

\newtheorem{lemma}[theorem]{Lemma}

\theoremstyle{remark}
\newtheorem*{remark}{Remark}

\everymath{\displaystyle}
\allowdisplaybreaks

\usepackage{authblk}

\makeatletter

\title{On the monoid of lexicographically minimal extensions}

\author[1]{Jonathan Caalim \thanks{jcaalim@math.upd.edu.ph} }

\author[2]{Yu-ichi Tanaka \thanks{tanaka\_yu-ichi@josogakuin.jp}}

\affil[1]{\small{University of the Philippines - Diliman, Quezon City, Philippines 1101}}
\affil[2]{\small{Joso Gakuin High School, Tsuchiura City, Ibaraki Prefecture, Japan 1010}}
\date{}

\begin{document}
\maketitle
\begin{abstract}
A sequence $(e_i)_{i \le m}$ of nonnegative integers $e_i$, where $m \in \mathbb{N}$ or $m =\infty$, is called a \textit{binomid index} if
\(\sum_{i=n-k+1}^{n} e_i\ge \sum_{i=1}^ke_i\)  for all $k, n \in \mathbb{N}$ such that $ 1\le k \le n < m$.
Infinite binomid indices give rise to binomid sequences (also known as Raney sequences) and generalized binomial coefficients. A finite binomid index $\eta$ can be extended to a unique lexicographically minimal infinite binomid index $\tilde{\eta}$. This lex-minimal extension $\tilde{\eta}$ is necessarily eventually periodic. In this research, we give a formula for the minimal period and provide an upper bound for the preperiod of $\tilde{\eta}$. We also show that the monoid of lex-minimal extensions is an inductive limit of finitely presented monoids.\\

\end{abstract}

\section{Introduction}
Let $\mathbb{N}:=\{1, 2, 3, \dots\}$, $\mathbb{N}_0:=\mathbb{N}\cup \{0\}$ and 
$(\mathbb{N}_0)^{\mathbb{N}}:=\{(e_1, e_2, \dots) \mid e_i \in \mathbb{N}_0\}$. 
A sequence $f=(f_n)_n=(f_1, f_2, \dots)$ of positive integers $f_n$ is a \textit{binomid sequence} (or Raney sequence) if, for all 
$n,k\in \mathbb{N}$ with $n \geq k$, the \textit{$f$-nomid coefficient}
\[ \begin{bmatrix} n \\ k \end{bmatrix}_f := \frac{f_n f_{n-1} \dots f_{n-k+1}}{f_kf_{k-1}\dots f_1}\]
is an integer (see \cite{shapiro2023, ando1999generalized}).
If $f$ is binomid, then there is a factorization of $f$ into prime factors $f(p) = (p^{\nu_p(f_n)})_n$, where $\nu_p(f_n)$ is the $p$-adic valuation of $f_n$, such that
\[f = \prod_{p \in \mathbb{P}} f(p)\]
with $\mathbb{P}$ denoting the set of prime numbers.
Consequently, each $f(p)$ is binomid and 
\[\begin{bmatrix} n \\ k \end{bmatrix}_f = \prod_{p \in \mathbb{P}} \begin{bmatrix} n \\ k \end{bmatrix}_{f(p)}.\]
Hence, it suffices to consider a sequence $f(p)=(p^{e_1}, p^{e_2}, \dots)$ 
for some $p\in \mathbb{P}$ and index  $(e_1, e_2, \dots)  \in (\mathbb{N}_0)^\mathbb{N}$.
An index $(e_i)_{i \le m}$, where $m \in \mathbb{N}$ or $m =\infty$, is called a \textit{binomid index} if,
  for all $k, n$ such that $ 1\le k \le n < m$,
 \[\sum_{i=n-k+1}^{n} e_i\ge \sum_{i=1}^ke_i.\]
We denote by $\mathbb{I}$ the set of all infinite binomid indices.

 A finite binomid index $\eta$ may always be extended into an infinite binomid index. Among the possible extensions in $\mathbb{I}$, we choose the unique lexicographically minimal one and denote this by $\tilde{\eta}$. We call $\tilde{\eta}$
 the \textit{lex-minimal extension} of $\eta$.

In \cite{caalimtanaka}, the authors showed that  the lex-minimal extension $\tilde{\eta}$ is eventually periodic. However, not all eventually periodic binomid indices are lex-minimal. Meanwhile, every periodic binomid index is lex-minimal. Moreover, the set $\mathbb{L}$ of all lex-minimal extensions is a commutative monoid under componentwise addition.

 For $l \in \mathbb{N}_0$,  let $\sigma^l$ be the (right) shift map that sends an infinite sequence
 $(e_1,e_2, \dots)$ to $( 0^l, e_1, e_2, \dots)$. From \cite{caalimtanaka}, we know that the elements of  the sum 
 \[\sum_{m \in \mathbb{N}_0} \sigma^m (\mathbb{L})\]
 are all eventually periodic binomid indices. We do not know whether the converse holds or not.

In this article, we reformulate the problem of computing the lex-minimal extension as an \textit{unbounded knapsack problem}.
We give a formula for the minimal period of $\tilde{\eta}$ and give an upper bound on the minimal preperiod.  Finally, we introduce a notion of atomic binomid indices to come up with a structure result (Theorem \ref{inductivelimit}) on the set of lex-minimal extensions. The atomic decomposition of a lex-minimal extension may provide a new perspective on solving an unbounded knapsack problem.

\section{Unbounded Knapsack Problem}

We begin this section by introducing some notations. 
Let $\mathbf{Sum}(\eta):=(S_i(\eta))_i=(S_i)_i$  be the sequence of partial sums of 
$\eta =(e_i)_i \in (\mathbb{N}_0)^{\mathbb{N}}$, i.e., $S_i = e_1 +e_2 + \dots +e_i$ for all $i \in \mathbb{N}$. 
Let $\mathbf{Ave}(\eta):=(A_i(\eta))_i=(A_i)_i$ be the sequence of averages
$A_i: = S_i/i$. For $k \in \mathbb{N}$, $\eta\big|_k$  denotes the truncated sequence $(e_i)_{i\le k}$. 


We have the following useful result that characterizes binomid indices and lex-minimal extensions.

\begin{proposition}\cite[Proposition 2.1]{caalimtanaka}\label{sumcondition} Let $\eta \in (\mathbb{N}_0)^\mathbb{N}$. 
Let $\Delta_\eta(k):=\min\{\delta_\eta(i,j)\mid i+j=k\}$ for $k\in \mathbb{N}$ with $k\ge 2$
 where 
$\delta_\eta (i,j):=S_{i+j} -S_i - S_j$ for $i, j \in \mathbb{N}$. Then
\begin{enumerate}
    \item $\eta$ is binomid if and only if $\delta_\eta (i,j) \ge 0$ for all $i, j \ge 1$;
    \item  $\eta$ is a lex-minimal extension of $\eta \big|_k$ 
if and only if $(\Delta_\eta(k'))_{k'> k}= 0^\infty$.
\end{enumerate}
\end{proposition}

 \begin{remark}
If $\eta \in \mathbb{I}$ and $k,n \in \mathbb{N}$, then $A_k \le A_{kn}$. Hence, $\lim_{i \to \infty} (A_i) $ exists (possibly as $\infty$). 
Moreover, $\delta$ is additive, i.e., for binomid indices $\eta$ and $\gamma$, \[\delta_{\eta+\gamma}(i, j) = \delta_\eta(i,j)+\delta_\gamma(i,j).\]
It follows that 
\[\Delta_{\eta+\gamma}(k) \ge \Delta_{\eta}(k)+\Delta_{\gamma}(k).\]
 \end{remark}

The process of computing the lex-minimal extension of a finite binomid index may be viewed as an \textit{unbounded knapsack problem} (UKP) (see \cite{Hu2009}). 

\begin{definition}
Let $v=(v_i)_{i\le m}$ and $w=(w_i)_{i\le m}$ be sequences in $(\mathbb{N}_0)^\mathbb{N}$. For $k \in \mathbb{N}$,  the unbounded knapsack problem $\mathbf{UKP}[v; w; k]$ asks to find a sequence $(n_i)_{i\le m}$ of  nonnegative integers that maximizes the sum 
   \[\sum_{i=1}^m n_i v_i\]
subject to the constraint
\[\sum_{i=1}^m n_i w_i \le k.\]
We denote the solution set to the problem also by $\mathbf{UKP}[v; w; k]$ and denote the maximal sum by $\phi(v; w; k)$. 
\end{definition}

Clearly, $\mathbf{UKP}[v; w; k]$ is nonempty if $k \ge \min \{w_i\}$.
If $w_i=1$ for some $i$, the condition $k \ge \sum_{i=1}^m n_i w_i$ can be replaced by $k = \sum_{i=1}^m n_i w_i$.

Given a $\mathbf{UKP}[v; w; k]$,
we can define the finite binomid index $\eta=(\eta_k)_{k\le \max{\{w_i\}}}$ where
$\eta(k) = \phi(v; w; k) -\phi(v; w; k-1)$.
Then $S(\tilde{\eta}) = (\phi(v; w; k))_{k\ge 1}$ for all  $k \in \mathbb{N}$.
Calculating the lex-minimal extension is equivalent to calculating maximal sums $\phi(v;w;k)$ arising from some unbounded knapsack problem. Indeed, we have the following result.

\begin{proposition}\label{cor:si_UKP}
 Let $\eta$ be a finite binomid index of length $m$ and let $\mathbf{Sum}(\tilde{\eta})=(S_i)_{i}$. 
 For all $k \ge 1$, 
there exists $(n_i)_{i\le m} \in \mathbf{UKP}[(S_i)_{i\le m}; (i)_{i\le m}; k]$ such that 
\[S_k = \sum_{i=1}^m n_i S_i\]
and \[k = \sum_{i=1}^m n_i i.\]
\end{proposition}

\begin{proof}
We prove the proposition by induction on $k \in \mathbb{N}$. If $k \le m$, the statement is clear. 
Suppose the result holds for all positive integers less than or equal to  $k$.
We show the result for $k+1$. Suppose $k+1 > m$. 
By Proposition \ref{sumcondition}, there are integers $k',k'' \ge 1$ such that  
$k+1 =k'+k''$ and 
$S_{k+1}=S_{k'}+S_{k''}$.
By the inductive hypothesis, there are nonnegative sequences $(n'_i)_{i\le m}$ and $(n''_i)_{i\le m}$ satisfying the statement of the result for $k'$ and $k''$. Let $n_i = n'_i+n''_i$ for $i\le m$. We  easily see that $(n_i)_{i\le m}$ satisfies the condition for $k+1$.
\end{proof}

\section{Preperiod and Period of $\tilde{\eta}$}

The lex-minimal extension $\tilde{\eta}$ of a finite binomid index $\eta$ is eventually periodic. 
Hence, we can write $\tilde{\eta}$ as $\eta_0 (\eta_1)^\infty$. We call $|\eta_1|$ the minimal period of $\tilde{\eta}$ if this length is minimal. Likewise, we call $|\eta_0|$ the minimal preperiod if this length is minimal. Of course, any integer greater than $|\eta_0|$ may be taken as a preperiod and any multiple $k |\eta_1|$, where $k \in \mathbb{N}$, may be taken as a period. 

In this section, we give an exact formula for the minimal period of $\tilde{\eta}.$
Throughout, we assume that $\eta$ is a finite binomid index of length $m$.
 We let $\mathbf{Sum}(\tilde{\eta})=(S_i)_i$ and $\mathbf{Ave}(\tilde{\eta})=(A_i)_{i}=(S_i/i)_{i}$. 
 We set $A_{\operatorname{max}} := \max \{A_i \mid 1\le i \le m\}$. 
 From \cite[Theorem 1]{caalimtanaka},  if $i \in \mathbb{N}$ such that $A_i =A_{\operatorname{max}}$, then  $i$ is a period of $\tilde{\eta}$.
 To proceed, we first prepare some lemmas.

\begin{lemma} \label{max_a} $A_{\operatorname{max}}=  \max \{A_i \mid i \in \mathbb{N}\}$
\end{lemma}

\begin{proof} 
By Proposition \ref{cor:si_UKP},  for $k>m$, there exists $(n_i)_{i\le m}\in \mathbf{UKP}[(S_i)_{i\le m}; (i)_{i\le m}; k]$ such that
\[k A_k = S_k = \sum_{i=1}^m n_i S_i=  \sum_{i=1}^m n_i  i   A_i\le \left(\sum_{i=1}^m n_i  i\right) A_{\operatorname{max}} = k A_{\operatorname{max}}.\]
Hence, $A_k \le A_{\operatorname{max}}$ for all $k\ge 1$. This implies that $A_{\operatorname{max}} = \max \{ A_k \mid k \in \mathbb{N}\}$.   
\end{proof}

\begin{lemma} \label{gcdequalsgcd}
$\mathbf{GCD}(\{i\le m \mid A_{\operatorname{max}}=A_i\}) = \mathbf{GCD}(\{i\in \mathbb{N} \mid A_{\operatorname{max}}=A_i\})$
\end{lemma}

\begin{proof}
Let $g = \mathbf{GCD}(\{i\le m \mid A_{\operatorname{max}}=A_i\})$. 
Let $k \in \mathbb{N}$ such that $A_k =A_{\operatorname{max}}$.
If $k \le m$, then $g\mid k$.
Suppose $k>m$. From the proof of Lemma \ref{max_a}, we have
\[\sum_{i=1}^m n_i  i   A_i= \sum_{i=1}^m n_i  i A_{\operatorname{max}} = kA_{\operatorname{max}}\]
since $kA_k=k A_{\operatorname{max}}$. 
Thus, 
if $A_i \neq A_{\operatorname{max}}$, then $n_i = 0$ necessarily because $A_{\operatorname{max}}$ is maximal.
It follows that
\[k =  \sum_{\substack{i=1 \\  A_i= A_{\operatorname{max}}}}^m n_i  i.\]
Hence, $g\mid k$. That is, $g=\mathbf{GCD}(\{i\in \mathbb{N}\mid A_{\operatorname{max}}=A_i\})$.
\end{proof}

\begin{lemma} \label{gcddividesperiod}
The minimal period of  $\tilde{\eta}$ is divisible by 
$\mathbf{GCD}(\{i\le m \mid A_{\operatorname{max}}=A_i\})$.
\end{lemma}

\begin{proof}
Since $\tilde{\eta}$ is a binomid index, it follows that $A_k \le A_{kn}$ for all $k,n \in \mathbb{N}$.
Therefore, if $k$ is in the set $ \{i\in \mathbb{N}\mid A_{\operatorname{max}}=A_i\}$, 
then, so is $nk$. 
Thus, we can choose $k\in \mathbb{N}$ sufficiently larger than both $m$ and the minimal preperiod of $\tilde{\eta}$
such that $A_k = A_{\operatorname{max}}$. 
Let $p$ be the minimal period of $\tilde{\eta}$ and
let $l:= k/\mathbf{GCD}(k,p)$. Since $k\mid lp$, we have
$A_{k+lp} = A_{\operatorname{max}}$.

Let $S$ be the sum of the $p$ entries of the periodic part of $\tilde{\eta}$. 
Then $S_{k+lp}=S_{k}+lS$. 
Also, 
\begin{eqnarray*}
    S_{k+lp}&=&(k+lp)A_{k+lp} \\
    &=&(k+lp)A_{\operatorname{max}} \\
    &=&k A_{\operatorname{max}} + lp A_{\operatorname{max}}\\
    &=& S_k+lp A_{\operatorname{max}}.
\end{eqnarray*}
Hence,  $lS=lp A_{\operatorname{max}}$, i.e., $S= p A_{\operatorname{max}}$.
It follows that
\begin{eqnarray*}
    (k+p)A_{k+p}&=&S_{k+p} \\
    &=&S_k+S\\
    &=&kA_{\operatorname{max}}+pA_{\operatorname{max}}\\
    &=&(k+p)A_{\operatorname{max}}.
\end{eqnarray*}
Thus, $A_{k+p} = A_{\operatorname{max}}$. 
So,  $k,  k+p \in \{i\in \mathbb{N}\mid A_{\operatorname{max}}=A_i\}$.
Therefore, $\mathbf{GCD}(\{i\in \mathbb{N}\mid A_{\operatorname{max}} = A_i\})$ divides $p= k+p-k$.
\end{proof}

\begin{theorem}\label{lengper}
 The minimal period of $\tilde{\eta}$ is $\mathbf{GCD}(\{i\le m \mid A_{\operatorname{max}}=A_i\})$. 
\end{theorem}

\begin{proof}
If $i, j \in \mathbb{N}$ are periods of $\tilde{\eta}$, then so is  $\mathbf{GCD}(i,j)$. 
Thus, $\mathbf{GCD}(\{i\le m \mid A_{\operatorname{max}}=A_i\})$ is a period of $\tilde{\eta}$. 
By Lemma \ref{gcddividesperiod}, it follows that $\mathbf{GCD}(\{i\le m \mid A_{\operatorname{max}}=A_i\})$ is the minimal period of $\tilde{\eta}$.
\end{proof}

\begin{theorem}\label{lengpreper}
If $k \in \{i \le m \mid A_i = A_{\operatorname{max}}\}$, then $km(m+1)/2$ is a preperiod of $\tilde{\eta}$.
\end{theorem}

\begin{proof} 
Let $k \le m$  such that $A_k=A_{\operatorname{max}}$. 
Let $q=km(m+1)/2$. Then $A_q = A_{\operatorname{max}}$ and $S_q = q A_{\operatorname{max}}$.
We show that 
\[S_{q+nk+r}=nS_k+S_{q+r}\]
for all $n \in \mathbb{N}$ and $0\le r < k$.
By Proposition \ref{cor:si_UKP},
there exists $(l_i)_{i\le m} \in \mathbf{UKP}[(S_i)_{i\le m}; (i)_{i\le m}; q+nk+r]$ such that
\[S_{q+nk+r}=\sum_{i=1}^m l_i S_{i}\]
and 
\[\sum_{i=1}^m il_i = q+nk+r.\]
For each $i \le m$,
there exist nonnegative integers $l'_i$ and $l''_i$ with $k|l'_i$ and $l''_i < k$ such that
$l_i = l'_i+l''_i$.
Then there is a natural number $t$ such that
\[\sum_{i=1}^m i l'_i = q+nk-tk\]
and 
\[\sum_{i=1}^m i l''_i= tk+r.\]
Since $l''_i < k$, 
the inequality $tk+r=\sum_{i=1}^m i l''_i  < \sum_{i=1}^m i k = km(m+1)/2 = q$.
Thus, $tk \le q.$
By Proposition \ref{cor:si_UKP},
\begin{eqnarray*}
S_{q+nk+r}&=&\sum_{i=1}^m l_i S_{i}\\
&=&\sum_{i=1}^m l'_i S_{i}+\sum_{i=1}^m l''_i S_{i} \\
&\le& S_{q+nk-tk}+\sum_{i=1}^m l''_i S_{i}\\
&\le& S_{q+nk+r}.
\end{eqnarray*}
In particular, $S_{q+nk+r}=S_{q+nk-tk}+\sum_{i=1}^m l''_i S_{i}$.
Meanwhile, $S_{q+nk-tk}=(q/k+n-t)S_k$ because $k\mid q+nk-tk$.
Then 
 \[S_{q+nk+r}=nS_k +(q/k-t)S_k+\sum_{i=1}^m l''_i S_{i}.\]
Define $(m_i)_{i\le m}$ where $m_k = l''_k + \frac{q}{k} - t$ and $m_i = l''_i$ for $i\neq k$. Then
\[S_{q+nk+r}=nS_{k}+\sum_{i=1}^m m_i S_i\]
and 
\[q+r = \sum_{i=1}^m im_i.\]
By Proposition \ref{cor:si_UKP},
\begin{eqnarray*}
 S_{q+nk+r}&=&nS_{k}+\sum_{i=1}^m m_i S_i \\
 &\le& nS_{k}+S_{q+r}\\
 &\le& S_{nk}+S_{q+r}\\
 &\le& S_{q+nk+r}.   
\end{eqnarray*}
Thus, $S_{q+nk+r}=nS_{k}+S_{q+r}$. This implies the periodicity of $(e_i)_{i>q}$. 
\end{proof}

The preperiod $km(m+1)/2$ in Theorem \ref{lengpreper} may not be optimal. Indeed, if $\tilde{\eta}$ is purely periodic, then the minimal preperiod is 0.


\section{Atomic Binomid Indices}

The set $\mathbb{L}$ of lex-minimal extensions is a monoid under componentwise addition. In this section, we provide a structure result for $\mathbb{L}$ by introducing the notion of atomic binomid indices. In particular, we show that the set of atomic lex-minimal extensions is the unique basis of the monoid $\mathbb{L}$.

We first prove the following result that asserts that if a lex-minimal extension is decomposable into binomid indices as a sum, then the summands, themselves, are lex-minimal extensions. 

\begin{proposition}\label{decomp2}
 Let $\eta$ be a lex-minimal extension. Let $\eta_1$ and $\eta_2$ be binomid indices such that $\eta = \eta_1 + \eta_2$. Then $\eta_1$ and $\eta_2$ are lex-minimal extensions, i.e., 
 \[(\mathbb{I}+\mathbb{I})\cap \mathbb{L} = \mathbb{L}+\mathbb{L}.\]
 Moreover, if $\eta$ is the lex-minimal extension of $\eta \big|_k$ for some $k \in \mathbb{N}$,  then $\eta_1$ and $\eta_2$ are  lex-minimal extensions of $\eta_1\big|_k$ and $\eta_2\big|_k$, respectively. 
\end{proposition}

\begin{proof}
Since $\eta_1$ and $\eta_2$ are binomid,
it follows that $\delta_{\eta_1}(i,j) \ge 0$ and $\delta_{\eta_2}(i,j) \ge 0$ for all $i,j \ge 1$. 
Suppose $\eta$ is the lex-minimal extension of $\eta \big|_k$ for some $k \in \mathbb{N}$. 
Then $\Delta_\eta(k')=0$ for all $k'>k$.
Moreover, for $k'>k$, there exist $i, j \ge 1$ such that $\delta_\eta(i,j)=0$ and $k'=i+j$. 
By additivity, $\delta_\eta = \delta_{\eta_1}+\delta_{\eta_2}$.
Hence, $\delta_{\eta_1}(i,j)=0$ and  $\delta_{\eta_2}(i,j)=0$. 
It follows that  $\Delta_{\eta_1}(k')=0$ and  $\Delta_{\eta_2}(k')=0$ for all $k' > k$. In other words, $\eta_1$ and $\eta_2$ are lex-minimal extensions of $\eta_1\big|_k$ and $\eta_2\big|_k$.
\end{proof}

\begin{definition}\label{atomicelement}
Let $\eta$ be a binomid index. We call $\eta$ atomic 
if there are no nonzero binomid indices  $\eta_1$ and $\eta_2$ such that $\eta = \eta_1+\eta_2$. 
We denote by $\mathbf{Atom}$ the set of all nonzero atomic binomid indices. 
\end{definition}

Let us consider some examples.

\begin{example}\label{exampleEm}
Let $m \in \mathbb{N}$. Clearly, $m^\infty$ is not atomic as
$m^\infty = 1^\infty+ \dots +1^{\infty}$ ($m$ summands of $1^\infty$).
Let $E_m := (0^{m-1},1)$.
It follows that $\widetilde{E_m}=(E_m)^\infty$ since $E_m$ is monotonically increasing (\cite[Proposition 3]{caalimtanaka}). Moreover, it is atomic. Indeed, if $E_m$ is a sum of two binomid indices, then one is necessarily zero. 

 Let $\eta=(0,1,1)$.
 Then $\widetilde{\eta}=(\eta)^\infty$. If $\widetilde{\eta}$ is not atomic,  then we decompose $\eta$ as
 $\eta= (0,0,1)+(0,1,0)$ necessarily. 
 Now, $\widetilde{(0,0,1)} = (0,0,1)^\infty$ and $\widetilde{(0,1,0)}=(0,1)^\infty$.
So, $\widetilde{(0,0,1)} + \widetilde{(0,1,0)}=(0,1,1,1,0,2)^\infty \neq \widetilde{\eta}.$ 
 Hence, $\tilde{\eta}$ is atomic. 
 Likewise, $(0^m,1,1)^\infty$ is an atomic lex-minimal extension for all  $m\in \mathbb{N}$. Thus, 
 \[(0^{m-1},1)^\infty,\ (0^m,1,1)^\infty \in \mathbf{Atom}\cap \mathbb{L}.\] 
\end{example}

\begin{example}\label{nonunique}
Let $\eta_1 = (0,1,1)$, $\eta_2=(0, 1, 2, 1, 1)$, $\eta_3=(0, 2, 1)$ and $\eta_4=(0, 0, 2, 0, 1)$.
Then
\begin{equation*}
\begin{aligned}
\widetilde{\eta_1}&=(0, 1, 1)^\infty\\
\widetilde{\eta_2}&=(0, 1, 2)(1)^\infty
\end{aligned}
\qquad
\begin{aligned}
\widetilde{\eta_3}&=(0,2)(1,1)^\infty\\
\widetilde{\eta_4}&=(0,0,2)(0,1,1)^\infty
\end{aligned}
\end{equation*}
\noindent
and each of these is atomic. Moreover, 
$\widetilde{\eta_1}+\widetilde{\eta_2}=\widetilde{\eta_3}+\widetilde{\eta_4}.$
Thus, the decomposition of a lex-minimal extension into a sum of atomic sequences is not unique.
\end{example}

To proceed, we first define two maps on the set of all finite binomid indices of a given length. 

\begin{definition}
  For \( k \in \mathbb{N} \), let \( \mathbb{I}_k \) be the set of all finite binomid indices of length \( k \). 
  Define the map \( L_k: \mathbb{I}_k \longrightarrow \mathbb{L} \) by $L_k(\eta) = \tilde{\eta}$.
  For $j,k \in \mathbb{N}$ with $j \le k$, 
  define the map \( L_{j,k}: \mathbb{I}_j \longrightarrow \mathbb{I}_k \)
  by
  \[L_{j,k}(\eta) = L_j(\eta)\big |_k = \tilde{\eta}\big |_k.\]
\end{definition}

\begin{remark}
As illustrated by Example \ref{exampleEm}, $L_{j,k}$ is not usually a homomorphism.
\end{remark}

Given a subset $A$ of a monoid $M$ with identity element 0, we denote by $\langle A \rangle$ the submonoid of $M$ given by the span 
\[\langle A \rangle :=\{\mu_1+\mu_2 + \dots + \mu_n \mid \mu_i \in A,\, n \in \mathbb{N}\} \cup \{ 0\}.\]
The subset $A$ is called a basis of $M$ if $M=\langle A \rangle$ and for all finite sequences \( (\beta_i)_{i\le n}\) over \( A \)  
  and \( (m_i)_{i \le n}\) over \(\mathbb{N}\cup \{0\} \), if
  \( \beta_j = \sum_{i=1}^n m_i \beta_i \) for some $1\le j \le n$, then \( m_j \ge 1 \).
Example \ref{nonunique} justifies this definition of the basis.

We adapt the definition of atomic binomid index in Definition \ref{atomicelement} to present the notion of an atomic element of a general monoid $M$, which is employed in the following useful result.

\begin{proposition}\label{uniquebasis}
 Let $M$ be a monoid and let $A$ be the set of all atomic elements of $M$. If every element of $M$ can be written as a finite sum of elements of $A$, then $A$ is the unique basis of $M$.
\end{proposition}

\begin{proof}
 We first demonstrate that $A$ is a basis of $M$.  
 By assumption, $A$ spans $M$. No element of $A$ can be decomposed into a nontrivial linear combination of other elements in \( A \) with nonnegative integer coefficients. Hence, \( A \) constitutes a basis for \( M \).
 
  Suppose $M$ has another basis \( B \) distinct from $A$. If there exists \( \nu \in A \setminus B \), then \( \nu \) must be expressible as a nontrivial linear combination of elements in \( B \). This, however, contradicts the atomic nature of \( \nu \). 
  Hence, \( A \subseteq B \). If there exists $\omega \in B \setminus A$, then $\omega$ can be written as $\omega=  \sum_{i=1}^n m_i \alpha_i$ where $\omega \neq \alpha_i \in A \subseteq B$ and $m_i\in \mathbb{N}_0$. Let $(\beta_i)_{i\le n+1}$ be a sequence in $B'$ such that $\beta_1 = \omega$ and $\beta_{i+1} = \alpha_i$ for all $i \le n$. Then $\beta_1 = 0 \cdot \beta_1 + \sum_{i=1}^{n} m_i\beta_{i+1}$. This contradicts the fact that $B$ is a basis. So, \( A = B \).
\end{proof}
\begin{proposition}
  The set $\mathbf{Atom}\cap \mathbb{L}$ is the unique basis of $\mathbb{L}$. 
  Moreover, $\mathbf{Atom}\cap  L_k(\mathbb{I}_k) $  is the unique basis of $\langle L_k(\mathbb{I}_k)\rangle $ for $k\in \mathbb{N}$.
\end{proposition}
\begin{proof}
 Let $\eta$ be an element of $\mathbb{L}$ such that $\eta \in L_k(\mathbb{I}_k)$ for some $k\in \mathbb{N}$.
 If $\eta=\sum_i \eta_i$ where $\eta_i \in \mathbb{L}$ for all $i$, then $\eta_i\in L_k(\mathbb{I}_k)$ by Proposition \ref{decomp2}. 
 Since the elements of $\mathbb{L}$ are sequences of nonnegative integers, the number of summands in the decomposition $\eta=\sum_i \eta_i$ is less than or equal to $S_k(\eta)$ and hence, the number of summands is finite.
 If we take the decomposition $\eta = \sum_i \eta_i$ with the maximal number of summands, then we see that all summands $\eta_i$ are in $\mathbf{Atom} \cap  L_k(\mathbb{I}_k) \subseteq \mathbf{Atom} \cap  \mathbb{L}$.
 Thus, the first statement follows from Proposition \ref{uniquebasis}.

To show the second statement, we prove that $\mathbf{Atom} \cap  L_k(\mathbb{I}_k) $ spans $\langle L_k(\mathbb{I}_k)\rangle$. 
Indeed, let $\eta\in \langle L_k(\mathbb{I}_k) \rangle$.
Then $\eta$ can be decomposed as a sum $\sum_i \eta_{i}$ where $\eta_{i}$ are elements of $L_k(\mathbb{I}_k)$.
By the above argument, each $\eta_{i}$ can be written as a finite sum of elements of $\mathbf{Atom} \cap  L_k(\mathbb{I}_k)$. 
This implies $\langle L_k(\mathbb{I}_k) \rangle$ satisfies thes the assumption of Proposition \ref{uniquebasis}.
Thus, $\mathbf{Atom} \cap  L_k(\mathbb{I}_k) $ is the unique basis of $\langle L_k(\mathbb{I}_k) \rangle$.
\end{proof}

In the latter part of this section, we
show that $\mathbf{Atom} \cap  L_k(\mathbb{I}_k) $ is finite. First, we mention the following remark.

\begin{remark}
For $k\in \mathbb{N}$, the space $\mathbb{I}_k$ is a polyhedral monoid defined by the conditions $\{S_{i+j}-S_i-S_j \ge 0 \mid i+j\le k\}$. 
 Thus,  $\mathbb{I}_k$ has a finite basis, which may 
be computed by using Integer Programming (see Bachem's method  in \cite{bachem1978theorem}).
\end{remark}

\begin{example}
 The set \( B= \{(1,1,1), (0,1,0), (0,0,1)\}\) is
 a basis of \( \mathbb{I}_3 \). 
  However, \( L_3(B) = \{(1)^\infty, (0,1)^\infty, (0,0,1)^\infty\} \) does not span $L_3(\mathbb{I}_3)$. 
  Indeed,  the lex-minimal extension \( \widetilde{(0,1,1)} = (0,1,1)^\infty\) cannot be written
  as a sum
  \(x (1)^\infty+ y(0,1)^\infty+z(0,0,1)^\infty\),
  where $x,y, z\in \mathbb{N}_0$.
\end{example}

 From the above example, we see that \( L_k(B) \) does not necessarily span
 \( L_k(\mathbb{I}_k) \) when $B$  is a basis of \( \mathbb{I}_k \). We have, instead, the following result.

\begin{lemma}\label{l(k)}
  For \( k  \in \mathbb{N}\),  there exists  \( l(k) \in \mathbb{N}\) with $l(k) \ge k$
  such that \( L_k(\mathbb{I}_k)\) is a subset of \(\langle L_{l(k)}(B)\rangle\) for any basis \( B \) of
  $\mathbb{I}_{l(k)}$.
\end{lemma}

\begin{proof}
  For \( \eta \in \mathbb{I}_k \), there exist finite sequences \( \eta_0 \) and \( \eta_1 \)  such that \( \tilde{\eta} = \eta_0(\eta_1)^\infty \). From Theorems \ref{lengper} and \ref{lengpreper}, if \( |\eta_0| \) and \( |\eta_1| \) are minimal, then they are bounded by \( k^3\) and \( k \), respectively. Moreover, if \( |\eta_1| \) is minimal, it divides \( \mathbf{LCM}(\{i\in \mathbb{N} \mid i \le k \})  \).  WLOG, we may assume that $ |\eta_0|=k^3$ and $|\eta_0\eta_1| = l(k):=k^3+\mathbf{LCM}(\{i\in \mathbb{N} \mid i \le k\})$.

    Let \( B \) be a basis of $\mathbb{I}_{l(k)}$. Then \( \tilde{\eta}\big|_{l(k)} \) can be expressed as a linear combination 
  \( \sum_{j} n_j \beta_j \) for some finite sequence \( (\beta_j)_j \) over \( B \) and a finite sequence  \( (n_j)_j \) of positive integers.  
  By Proposition \ref{sumcondition}, 
  we have \( \Delta_{\beta_j}(m) \le \Delta_{\tilde{\eta}}(m) \) for all \( j \) and \( 2\le m \le l(k) \). Therefore, \( \Delta_{\beta_j}(m) = 0 \) for all \( l(k)\ge m > k \) as 
  \( \Delta_{\tilde{\eta}}(m) = 0 \) for these values of \( m \).

Write \(\widetilde{\beta_j\big|_k} =\widetilde{\beta_j} = \beta_{(j,0)}(\beta_{(j,1)})^\infty \) 
  with $|\beta_{(j,0)}|=k^3$ and $|\beta_{(j,0)}\beta_{(j,1)}| = l(k)$.
  Then \( \beta_j = \beta_{(j,0)}\beta_{(j,1)} \).
   Since \( \tilde{\eta}\big|_{l(k)} = \eta_0\eta_1 \), we have \( \eta_0 = \sum_{j} n_j \beta_{(j,0)} \) and \( \eta_1 = \sum_{j} n_j \beta_{(j,1)} \). Thus, 
  \[
  \tilde{\eta} = \eta_0(\eta_1)^\infty = \sum_j n_j \beta_{(j,0)}(\beta_{(j,1)})^\infty = \sum_j n_j \widetilde{\beta_j}.
  \]
  Therefore, \( L_k(\mathbb{I}_k) \) is a subset of the set spanned by \( L_{l(k)}(B) \). 
\end{proof}

In what follows, we show that that the monoid $\langle L_k(\mathbb{I}_k)\rangle$ is finitely presented. 
Recall that, by Lemma \ref{l(k)}, there exists $l(k) \in \mathbb{N}$ with $l(k)\geq k$ such that $ \langle L_k(\mathbb{I}_k)\rangle$ is a submonoid of  
$\langle L_{l(k)}(B)\rangle$  for any basis $B$ of $\mathbb{I}_{l(k)}$. For $k \in \mathbb{N}$, set $(k,l(k)]:=\{i \in \mathbb{N} \mid k < i \le l(k)\}$. Let $\mathbf{Path}$ be the set of all maps 
 \[\pi:(k,l(k)] \longrightarrow \{i \in \mathbb{N} \mid i < l(k)\}\]
 such that $\pi(i) < i$ for all $i  \in (k,l(k)]$.
By Proposition \ref{cor:si_UKP}, we have the following lemma.

\begin{lemma}  If $\eta \in L_{k,l(k)}(\mathbb{I}_k)$, 
then there is a map $\pi \in \mathbf{Path}$ such that 
\[S_{i}-S_{\pi(i)}-S_{i-\pi(i)} = 0\]
for all $ i \in (k,l(k)]$.
\end{lemma}

\begin{lemma}\label{path}
The monoid $\langle L_k(\mathbb{I}_k)\rangle$ is finitely generated for all $k \in \mathbb{N}$. In particular,  its basis $\mathbf{Atom} \cap  L_k(\mathbb{I}_k) $ is finite.
\end{lemma}

\begin{proof}
Let $l(k) \in \mathbb{N}$ such that $l(k)\geq k$ and 
$ \langle L_k(\mathbb{I}_k)\rangle$ is a submonoid of  
$\langle L_{l(k)}(B)\rangle$  for any basis $B$ of $\mathbb{I}_{l(k)}$.
Consider the restriction map $\operatorname{Res}_{l(k)}:\langle L_{l(k)}(B)\rangle \longrightarrow \mathbb{I}_{l(k)}$ given by $\operatorname{Res}_{l(k)}(\eta)=\eta |_{l(k)}$.
Let $B=\{\beta_i \mid 1\le i \le n\}$. 
If $\eta \in \langle L_{l(k)}(B)\rangle$, then there exist $m_i \in \mathbb{N}\cup \{0\}$  such that
$\eta = \sum_{i=1}^n m_i L_{l(k)}(\beta_i)$. Observe that 
\begin{eqnarray*}
\operatorname{Res}_{l(k)}(\eta)&=&\sum_{i=1}^n m_i \operatorname{Res}_{l(k)}(L_{l(k)}(\beta_i))\\
             &=&\sum_{i=1}^n m_i \beta_i \in \mathbb{I}_{l(k)}. 
\end{eqnarray*}
In other words, $\operatorname{Res}_{l(k)}$  maps $\sum_{i=1}^n m_i L_{l(k)}(\beta_i)$  to $\sum_{i=1}^n m_i \beta_i$. Thus, $\operatorname{Res}_{l(k)}$ is a monoid isomorphism between $\langle L_{l(k)}(B)\rangle$ and $\mathbb{I}_{l(k)}$.

Next, we show that $\operatorname{Res}_{l(k)}(\langle L_k(\mathbb{I}_k)\rangle)\cong \langle L_k(\mathbb{I}_k)\rangle$. 
By Lemma \ref{l(k)}, $\langle L_k(\mathbb{I}_k)\rangle \subseteq \langle L_{l(k)}(B)\rangle$. 
Restricting the domain of the map $\operatorname{Res}_{l(k)}$ to the submonoid $\langle L_k(\mathbb{I}_k)\rangle$, we see that $\langle L_k(\mathbb{I}_k)\rangle$ is isomorphic to $\operatorname{Res}_{l(k)}(\langle L_k(\mathbb{I}_k)\rangle)$ as monoids.

 Let $\eta \in L_{k,l(k)}(\mathbb{I}_k)= \operatorname{Res}_{l(k)}(L_k(\mathbb{I}_k))$. By Lemma \ref{path}, there exists a map $\pi \in \mathbf{Path}$ such that, for all $ i \in (k,l(k)]$, 
 \[S_{i}-S_{\pi(i)}-S_{i-\pi(i)} = 0.\]
Moreover, if $\eta \in \mathbb{I}_{l(k)}$ satisfies the equation 
$ S_{i}-S_{\pi(i)}-S_{i-\pi(i)} = 0$ for all $i \in (k,l(k)]$ for some $\pi \in \mathbf{Path} $, 
then there exists $\mu \in \mathbb{I}_{k}$ such that 
\[\eta = L_{k,l(k)}(\mu).\]
 Thus, for $\pi\in \mathbf{Path}$, the set 
 \[\operatorname{Poly}(\pi):=\{\eta \in \mathbb{I}_{l(k)}\mid S_{i}-S_{\pi(i)}-S_{i-\pi(i)} = 0 \text{ for all }i\in (k, l(k)]\},\]
 which is clearly closed under addition, is 
 a submonoid of $\operatorname{Res}_{l(k)}(\langle L_k(\mathbb{I}_k)\rangle)$. 
 Since the condition $S_{i}-S_{\pi(i)}-S_{i-\pi(i)}=0$ is equivalent to the following system of inequalities
 \[\begin{cases}
 \phantom{-(}S_{i}-S_{\pi(i)}-S_{i-\pi(i)}\phantom{)}  \geq 0  \\
 -(S_{i}-S_{\pi(i)}-S_{i-\pi(i)})\geq 0, 
 \end{cases}\]
 the set $\operatorname{Poly}(\pi)$ is a polyhedral submonoid of $\mathbb{I}_{l(k)}\subset \mathbb{N}^{l(k)}$.
 By a generalized Minkowski's theorem of Hilbert (see \cite{hilbert1890ueber, bachem1978theorem, jeroslow1978some}), $\operatorname{Poly}(\pi)$ is finitely generated.
 Likewise, $\sum_{\pi \in \mathbf{Path}} \operatorname{Poly}(\pi)$ is finitely generated because $\mathbf{Path}$ is finite.
 Recall that $\operatorname{Res}_{l(k)}(L_k(\mathbb{I}_k)) \subset \bigcup_{\pi \in \mathbf{Path}} \operatorname{Poly}(\pi)$. Hence, we have 
 \[\sum_{\pi \in \mathbf{Path}} \operatorname{Poly}(\pi) =\operatorname{Res}(\langle L_k(\mathbb{I}_k)\rangle).\]
 It follows that $\operatorname{Res}_{l(k)}(\langle L_k(\mathbb{I}_k)\rangle)$ is finitely generated. 
 Thus, $ \langle L_k(\mathbb{I}_k)\rangle$ is finitely generated.
\end{proof}

R\'edei \cite{redei2014theory} proved that when a commutative monoid is finitely generated, it is also finitely presented. Hence, we have the following result.

\begin{lemma}
The monoid $\langle L_k(\mathbb{I}_k)\rangle$ of $\mathbb{L}$ is finitely presented for all $k \in \mathbb{N}$. 
\end{lemma}

Finally, we have the following result.
 
\begin{theorem}\label{inductivelimit}
The set $\mathbb{L}$ of lex-minimal extensions is an inductive limit of the increasing system $(\langle L_k(\mathbb{I}_k)\rangle)_{k \in \mathbb{N}}$ of finitely presented monoids. 
Moreover,  $(\mathbf{Atom}\cap  L_k(\mathbb{I}_k))_{k \in \mathbb{N}}$ is an increasing system of generators of $(\langle L_k(\mathbb{I}_k)\rangle)_{k \in \mathbb{N}}$.
\end{theorem}


\begin{thebibliography}{9}

\bibitem{bachem1978theorem}Bachem, A. The theorem of Minkowski for polyhedral monoids and aggregated linear diophantine systems. {\em Optimization and Operations Research: Proceedings of a Workshop Held at the University Of Bonn, October 2–8, 1977}. pp. 1-13 (1978)
\bibitem{caalimtanaka}Caalim, J. \& Tanaka, Y. Lexicographically Minimal Extension of a Finite Binomid Index. {\em Integers}. \textbf{24 \#A53}, 1-13 (2024)
\bibitem{ando1999generalized}Ando, S. \& Sato, D. On the generalized binomial coefficients defined by strong divisibility sequences. {\em Applications of Fibonacci Numbers: Volume 8}. pp. 1-10 (1999)
\bibitem{facchini2002direct}Facchini, A. Direct sum decompositions of modules, semilocal endomorphism rings, and Krull monoids. {\em Journal of Algebra}. \textbf{256}, 280-307 (2002)
\bibitem{jeroslow1978some}Jeroslow, R. Some basis theorems for integral monoids. {\em Mathematics of Operations Research}. \textbf{3}, 145-154 (1978)
\bibitem{Hu2009}Hu, T., Landa, L. \& Shing, M. The unbounded knapsack problem. {\em Research Trends in Combinatorial Optimization: Bonn 2008}. pp. 201-217 (2009)
\bibitem{hilbert1890ueber}Hilbert, D. Ueber die Theorie der algebraischen Formen. {\em Mathematische Annalen}. \textbf{36}, 473-534 (1890)
\bibitem{shapiro2023}Shapiro, D. Divisibility properties for integer sequences. {\em Integers}. \textbf{23 \#A57}, 1-20 (2023)
\bibitem{redei2014theory}Rédei, L. The theory of finitely generated commutative semigroups. (Elsevier,2014)

\end{thebibliography}
\end{document}